\documentclass[12pt]{amsart}
\usepackage{amsfonts}
\usepackage{amssymb}
\usepackage[bulgarian, english]{babel}
\newtheorem{thm}{Theorem}[section]
\newtheorem{cor}[thm]{Corollary}
\newtheorem{lem}[thm]{Lemma}

\theoremstyle{definition}
\newtheorem{defn}[thm]{Definition}
\theoremstyle{remark}

\newcommand{\X}{\mathfrak X}

\begin{document}

\title[Holomorphic hypersurfaces of K\"ahler manifolds with Norden metric]
{Holomorphic hypersurfaces of K\"ahler manifolds with Norden metric}%
\author{Georgi Ganchev, Kostadin Gribachev and Vesselka Mihova}%
\address{Bulgarian Academy of Sciences, Institute of Mathematics and Informatics,
Acad. G. Bonchev Str. bl. 8, 1113 Sofia, Bulgaria}%
\email{ganchev@math.bas.bg}%
\address{} %
\address{Faculty of Mathematics and Informatics, University of Sofia,
J. Bouchier Str. 5, (1164) Sofia, Bulgaria}
\email{mihova@fmi.uni-sofia.bg}

Universite de Plovdiv ``Paissi Hilendarski",

Travaux scientifiques, Mathematiques

vol. 23, no. 2 (1985), 221 - 237 \vskip 10mm

\subjclass{Primary 53B05}%
\keywords{K\"ahlerian manifolds with Norden metric, holomorphic hypersurfaces,
totally real sectional curvatures, classification theorem.}%
\maketitle

\thispagestyle{empty}
\begin{abstract}
We study K\"ahlerian manifolds with Norden metric $g$ and develop the theory of
their holomorphic hypersurfaces with constant totally real sectional curvatures.
We prove a classification theorem for the holomorphic hypersurfaces of $(\mathbb{R}^{2n+2},
g, J)$ with constant totally real sectional curvatures.
\end{abstract}

\vskip 6mm

On an almost Hermitian manifold $(M, h, J)$ the restriction of the almost complex structure
$J$ on the tangent space at an arbitrary point of $M$ is an isometry with respect to the metric
$h$. This relation between both structures $J$ and $h$ determines the almost Hermitian geometry.
On an almost complex manifold $(M, g, J)$ with Norden metric $g$ the restriction of the almost
complex structure $J$ on the tangent space at an arbitrary point of $M$ is an anti-isometry with
respect to the metric $g$. This relation between both structures $J$ and $g$ determines the
geometry of the almost complex manifolds with Norden metric. The basic classes of such manifolds
have been given in \cite {2}. The most important of these classes is the class of K\"ahler
manifolds with Norden metric. A four-dimensional example of a K\"ahler manifold with Norden
metric has been given in \cite {4}. Another approach to K\"ahler manifolds with Norden metric
has been used in \cite {5}. In \cite {3}, there has been proved that the four-dimensional sphere
of Kotel'nikov-Study carries a structure of a K\"ahler manifold with Norden metric.

An essential problem in the differential geometry of K\"ahler manifolds is the investigation of
the complex space forms (K\"ahler manifolds of constant holomorphic sectional curvatures).
The corresponding problem in the differential geometry of K\"ahler manifolds with Norden metric
is to make a study of the manifolds with constant totally real sectional curvatures $\nu$ and
$\tilde{\nu}$. In this paper we give examples of K\"ahler manifolds with Norden metric having
prescribed curvatures $\nu$ and $\tilde\nu$. We prove a classification theorem for the
holomorphic hypersurfaces of $(\mathbb{R}^{2n+2}, g, J)$ with constant totally real sectional
curvatures $\nu$ and $\tilde\nu$.
\vskip 2mm

\section{Algebraic preliminaries}
Let $g_0$ be the usual scalar product on $\mathbb{R}^n$, i.e.
$$g_0(\xi, \eta)=\delta_{ij}\, \xi^i \,\eta^j,\quad \xi=(\xi^1,...,\xi^n)\in\mathbb{R}^n,\;
\eta=(\eta^1,...,\eta^n)\in\mathbb{R}^n.$$

The complex linearization
$G_0$ of $g_0$ is a symmetric complex bilinear form on
$\mathbb{C}^n$. The usual identification $r$ of $\mathbb{C}^n$ with
$\mathbb{R}^{2n}$ is given by
$$r:\;z=(z^1,...,z^n)\in\mathbb{C}^n\,\mapsto\, r(z)=Z=(x^1,...,x^n; y^1...,y^n)\in
\mathbb{R}^{2n},$$
where $z^k=x^k+i\,y^k\; (k=1,..., n)$. The canonical complex structure $J$ on $\mathbb{R}^{2n}$
(induced by the complex structure on $\mathbb{C}^n$) is determined by the matrix
$$\left(\begin{array}{cc}
0& I_n\\
[2mm]
-I_n & 0
\end{array}\right)$$
with respect to the natural basis of $\mathbb{R}^{2n}$. $G_0$ induces two metrics on
$\mathbb{R}^{2n}$:
$$g:={\rm Re}\,G_0,\quad \tilde g:={\rm Im}\,G_0.$$

If $Z=(x^1,...,x^n; y^1...,y^n)\in \mathbb{R}^{2n}$, then
$$\begin{array}{ll}
g(Z,Z)=&(x^1)^2+...+(x^n)^2-(y^1)^2-...-(y^n)^2,\\
[2mm]
\tilde g(Z,Z)=& 2(x^1\,y^1+...+x^n\,y^n).
\end{array}\leqno (1.1)$$

For all $Z,\, W$ in $\mathbb{R}^{2n}$ the metric $g$ and the complex structure $J$ on
$\mathbb{R}^{2n}$ are related by the equality $g(JZ, JW)=-g(Z, W)$. The metric $\tilde g$ is
said to be an associated (to $g$) metric because of $\tilde g(Z, W)=g(JZ, W)$. Both metrics
are of signature $(n, n)$.

The subgroup $O(n; \mathbb{C})$ of ${\rm GL}(n, \mathbb{C})$ preserving $G_0$ consists of
the matrices $\alpha$ such that $\alpha\,\alpha^t = I_n$. If $r$ is the real representation of
${\rm GL}(n, \mathbb{C})$, then $r(O(n, \mathbb{C}))=r({\rm GL}(n, \mathbb{C}))\cap O(n, n)$,
i.e. for $\alpha= A+i\,B\in O(n, \mathbb{C})$
$$r(\alpha)=\left(\begin{array}{rc}
A & B\\
[2mm]
-B & A
\end{array}\right),\quad A\,A^t- B\,B^t = I_n,\quad A\,B^t+B\,A^t=0.$$

Thus, the group $r(O(n, \mathbb{C}))$ consists of the matrices of ${\rm GL}(2n, \mathbb{R})$
preserving both structures $J$ and $g$ on $\mathbb{R}^{2n}$.

Now, let $(V, g, J)$ be a $2n$-dimensional real vector space with metric $g$ and complex structure
$J$, so that $g(Jx, Jy)=-g(x, y)$ for all $x, y$ in $V$. The associated metric $\tilde g$ is
given by $\tilde g(x, y)=g(J x, y),\; x, y \in V$. The group $r(O(n, \mathbb{C}))$ is the
group preserving both structures $g$ and $J$ on $V$. There exist bases of $V$ of the type
$\{x_1, ..., x_n; Jx_1, ..., Jx_n\}$ such that $g(x_i, x_j)=-g(Jx_i, Jx_j)= \delta_{ij},\;
g(x_i, Jx_j)=0;\; i, j=1, ..., n$. Such bases are called adapted bases with respect to $g$.
The linear transformation of an adapted basis to another one is given by a matrix of
$r(O(n, \mathbb{C}))$.

\begin{defn}
The matrix $S\in {\rm GL} (2n, \mathbb{R})$ is said to be holomorphic symmetric (h-symmetric) if $S\circ J= J\circ S$
and $g(S x, y)=g(x, S y)$ for all $x, y\in V$.
\end{defn}

It follows immediately that $S\in {\rm GL} (2n, \mathbb{R})$ is h-symmetric if and only if
$$g(S x, y)= g(x, Sy),\quad \tilde g(S x, y)= \tilde g(x, S y),\; x, y \in V.$$
\begin{defn}
Let $S\in {\rm GL} (2n, \mathbb{R})$ and $S\circ J= J\circ S$. The vector $x \in V$ is said to be
holomorphic proper (h-proper) for $S$ if
$$S x=\lambda x+ \mu J x;\quad \lambda,\, \mu \in \mathbb{R}.$$
\end{defn}

It is clear that $S \in {\rm GL} (2n, \mathbb{R})$ is h-symmetric if and only if
$r^{-1}(S)$ is symmetric. This implies
\begin{lem}\label {L: 1.1}
Let $S \in {\rm GL} (2n, \mathbb{R})$ be h-symmetric. Then there exists in $V$ an adapted basis
of h-proper vectors of $S$. With respect to such a basis $S$ has the form
$$S=\left(\begin{array}{rc}
A & B\\
[2mm]
-B & A
\end{array}\right),\quad
A=\left(\begin{array}{ccc}
\lambda_1 && 0\\
&\ddots&\\
0 && \lambda_n
\end{array}\right),\quad
B=\left(\begin{array}{ccc}
\mu_1 && 0\\
&\ddots&\\
0 && \mu_n
\end{array}\right),$$
where $\lambda_i,\, \mu_i \in \mathbb{R},\; i=1,...,n$.
\end{lem}
\vskip 2mm

\section{Holomorphic umbilical hypersurfaces of K\"ahler manifolds with Norden metric}

Let $(M, J)\; (\dim M= 2n)$ be an almost complex manifold. A metric $g$ on $M$ is said to be a
Norden metric if $g(J x, J y)=- g(x, y)$ for all vector fields $x, y$ on $M$. The metric $g$
is necessarily of signature $(n, n)$. An almost complex manifold with Norden metric
$(M, g, J)$ is said to be a K\"ahler manifold with Norden metric if $\nabla J =0$ with respect
to the Levi-Civita connection of $g$ \cite {2}. The tangent space $T_p M,\; p\in M$, of an almost
complex manifold with Norden metric is a vector space with Norden metric $g$ and complex
structure $J$ as considered in section 1.

Let $(M, g, J),\; (\dim M=2n)$ be a K\"ahler manifold with Norden metric. The algebra of the
differentiable vector fields on $M$ is denoted by $\X M$. The curvature tensor $R$ of $g$ is
given by $R(X, Y) Z= \nabla_X \nabla _Y Z-\nabla_Y \nabla _X Z-\nabla_{[X, Y]} Z$ for all
$X, Y, Z \in \X M$. The curvature tensor $R$ of type $(0, 4)$ is defined by $R(x, y, z, u)=
g(R(x, y) z, u)$ for all $x, y, z, u$ in $T_p M,\; p \in M$, and has the property
$$R(x,y,z,u)=-R(x,y,J z, Ju);\quad x,y,z,u \in T_p M,\; p\in M.$$

This implies that the tensor $\tilde R$ defined by $\tilde R(x, y, z, u)=R(x, y, z, J u)$
has the property $\tilde R(x, y, z, u)=-\tilde R(x, y, u, z).$

In the theory of K\"ahler manifolds with Norden metric the following tensors are essential:
$$\begin{array}{l}
\pi_1(x, y, z, u)= g(y, z)\, g(x, u)- g(x, z)\, g(y, u),\\
[2mm]
\pi_2(x, y, z, u)= \tilde g(y, z)\,\tilde g(x, u)-\tilde g(x, z)\,\tilde g(y, u),\\
[2mm]
\pi_3(x, y, z, u)= -g(y, z)\,\tilde g(x, u)+ g(x, z)\,\tilde g(y, u)
-\tilde g(y, z)\, g(x, u)+\tilde g(x, z)\, g(y, u),
\end{array}$$
$x,y,z,u \in T_p M,\; p\in M$.

Every non degenerate with respect to $g$ 2-plane $\beta$ in $T_p M,\; p \in M$, has two
sectional curvatures
$$K(\beta; p)=\frac{R(x,y,y,x)}{\pi_1(x,y,y,x)},\quad
\tilde K(\beta; p)=\frac{\tilde R(x,y,y,x)}{\pi_1(x,y,y,x)},$$
where $\{x, y\}$ is a basis of $\beta$.

A 2-plane $\beta$ in $T_p M,\; p \in M$ is said to be totally real if $\beta \neq J \beta$
and $ \beta \perp J \beta$ with respect to $g$.

\begin{thm} \cite {1} \label {T: 2.1}
Let $M$ \;$(\dim M= 2n \geq 4)$ be a K\"ahler manifold with Norden metric. $M$ is of
pointwise constant totally real sectional curvatures $\nu$ and $\tilde \nu$, i.e.
$$K(\beta; p)= \nu(p),\quad \tilde K(\beta; p)=\tilde \nu(p)$$
for an arbitrary non degenerate totally real 2-plane $\beta$ in $T_p M,\; p\in M$,
if and only if
$$R=\nu\, (\pi_1-\pi_2)+\tilde \nu \,\pi_3.$$

Both functions  $\nu$ and $\tilde \nu$ are constant if $M$ is connected and $\dim M\geq 6$.
\end{thm}

Let $(M', g, J)\; (\dim M'=2n+2)$ be a K\"ahler manifold with Norden metric and let $\nabla'$,
$R'$ be the Levi-Civita connection and its curvature tensor, respectively. A submanifold
$M\; (\dim M=2n)$ is said to be a holomorphic hypersurface of $M'$ if the restriction of
$g$ on $M$ has a maximal rank and $J T_p M = T_p M,\; p\in M$. We denote the restrictions of
$g$ and $J$ on $M$ by the same letters. Then, $(M, g, J)$ is an almost complex manifold with
Norden metric. There exist locally vector fields $\xi$ and $J \xi$ normal to $M$, such that
$$g(\xi, \xi)=-g(J\xi, J\xi)=1,\quad g(\xi, J\xi)=0.\leqno (2.1)$$

In fact, we can choose locally vector fields $\eta$ and $J\eta$ normal to $M$ such that
$g(\eta, \eta)= -g(J\eta, J\eta)=1$. Putting $g(\eta, J\eta)=\sinh t$, we obtain that
$$\xi=\frac{1}{\cosh t}\left\{\left(\cosh{\frac{t}{2}}\right)\eta+
\left(\sinh{\frac{t}{2}}\right)J\eta\right\},\;
J\xi=\frac{1}{\cosh t}\left\{\left(-\sinh{\frac{t}{2}}\right)\eta+
\left(\cosh{\frac{t}{2}}\right)J\eta\right\}$$
satisfy (2.1).

If $\nabla$ is the induced Levi-Civita connection on $M$, then the Gau\ss\; and Weingarten
formulas are
$$\begin{array}{lll}
\nabla'_X Y&= \nabla_X Y + \sigma(X, Y);& X, Y \in \X M,\\
[2mm]
\nabla'_X \xi &=  - A X + D_X \xi;& X \in \X M,
\end{array}$$
where $\sigma$ is the second fundamental form on $M$, $A=A_{\xi}$ is the second
fundamental tensor with respect to $\xi$, and $D$ is the normal connection on $M$.
Because of $\nabla' g=\nabla' J=0$ and
$$g(\sigma(x, y), \xi)= g(Ax, y)= g(x, Ay);\quad x, y \in T_p M,\; p\in M,$$
we obtain
$$\begin{array}{rll}
\sigma(x, y)&=& g(A x, y)\,\xi - \tilde g(A x, y)\,J \xi;\\
[2mm]
\sigma(x, J y)&=&\sigma(J x, y)= J\sigma(x, y);\\
[2mm]
A_{J\xi}&=& A\circ J = J\circ A;\\
[2mm]
D_x \xi &=& D_x J\xi = 0,\quad \left(\nabla_x J\right)\,y = 0
\end{array}\leqno (2.2)$$
for arbitrary $x, y \in T_p M,\; p\in M$.

Thus, every holomorphic hypersurface $M$ of $M'$ is also a K\"ahler
manifold with Norden metric and the Gau\ss\; and Weingarten formulas for $M$ are
$$\begin{array}{lll}
\nabla'_X Y&= \nabla_X Y + g(A X, Y)\,\xi - \tilde g(A X, Y)\,J \xi;& X, Y \in \X M,\\
[2mm]
\nabla'_X \xi &=  - A X;& X \in \X M.
\end{array}\leqno (2.3)$$

From now on, $(M', g, J)\; (\dim M' = 2n+2)$ will stay for a K\"ahler manifold with Norden
metric and $(M, g, J)$ will stay for a holomorphic hypersurface of $M'$.

\begin{lem}\label {L: 2.1}
Let $R'$ and $R$ be the curvature tensors of $M'$ and $M$ respectively. Then
$$\begin{array}{ll}
R'(x, y, z, u)&= R(x, y, z, u)-\pi_1(Ax, Ay, z, u)+\pi_2(Ax, Ay, z, u);\\
[2mm]
(R'(x, y)z)^{\perp}&=g((\nabla_x A) y-(\nabla_y A) x, z)\,\xi-
\tilde g((\nabla_x A) y-(\nabla_y A) x, z)\,J\xi;\\
[2mm]
R'(x, y)\xi&=(\nabla_y A) x-(\nabla_x A) y
\end{array}\leqno (2.4)$$
for all $x, y, z, u \in T_p M,\; p \in M$.
\end{lem}

From Lemma \ref {L: 2.1} in a straightforward way we obtain
\begin{lem}\label {L: 2.2}
Let $\rho'$ and $\rho$ be the Ricci tensors of $M'$ and $M$ respectively. Then
$$\begin{array}{ll}
\rho'(x, y)=&\rho(x, y)- ({\rm trace}\,A)\,g(Ax, y)+({\rm trace}\,(A\circ J))\,g(Ax, Jy)\\
[2mm]
&+ 2 g(A^2 x, y)+2 R'(\xi, x, y, \xi).
\end{array}\leqno (2.5)$$
\end{lem}

Further, taking into account Theorem \ref {T: 2.1} and Lemma \ref {L: 2.2}, we get
\begin{lem}\label {L: 2.3}
Let $M'$ be of constant totally real sectional curvatures $\nu\,'$ and ${\tilde \nu}\,'$. Then
$$\begin{array}{ll}
(\nabla_x A) y &= \nabla_y A) x;\\
[2mm]
\rho(x, y)&= 2(n-1)(\nu\,'\,g(x, y)- \tilde\nu\,'\,g(x, Jy))\\
[2mm]
&+({\rm trace}\,A)\,g(Ax, y)-({\rm trace}\,(A\circ J))\,g(Ax, Jy)-2 g(A^2 x, y)
\end{array}\leqno (2.6)$$
for arbitrary vectors $x, y \in T_p M,\; p \in M$.
\end{lem}

Let $H$ be the mean curvature vector  on $M$, i.e. $H=\displaystyle{\frac{{\rm trace\,\sigma}}
{2n}}$. Taking into account (2.2), we get
$$\begin{array}{rll}
H&=\displaystyle{\frac{{\rm trace\,\sigma}}{2n}}&
\displaystyle{=\frac{1}{2n}\,(({\rm trace\,A})\,\xi
-({\rm trace\,(A\circ J)})\,J\xi)};\\
[2mm]
JH&=\displaystyle{\frac{{\rm trace\,(J\circ\sigma)}}{2n}}&
\displaystyle{= \frac{1}{2n}\,(({\rm trace\,(A\circ J)})\,\xi
+({\rm trace\,A})\,J\xi)};
\end{array}\leqno (2.7)$$
$$\begin{array}{ll}
g(H, H)&=\displaystyle{\left(\frac{{\rm trace\,A}}{2n}\right)^2-
\left(\frac{{\rm trace\,(A\circ J)}}{2n}\right)^2,}\\
[3mm]
\tilde g(H, H)&=\displaystyle{2\;\frac{{\rm trace\,A}}{2n}\;\frac{{\rm trace\,(A\circ J)}}{2n}.}
\end{array}\leqno (2.8)$$

These formulas imply that $M$ is minimal if and only if
${\rm trace\,A}={\rm trace\,(A\circ J)}=0$.

\begin{defn}
The manifold $M$ is said to be holomorphically umbilical (h-umbilical) if at every point of $M$
$$\sigma=g\,\frac{{\rm trace}\,\sigma}{2n}-\tilde g \,\frac{{\rm trace}\,(J\circ\sigma)}{2n}
=g\,H- \tilde g\,JH.
\leqno (2.9)$$
\end{defn}

It is clear that if $M$ is h-umbilical and minimal, then $M$ is totally geodesic.

Let $\eta$ be an arbitrary vector, normal to $M$ at $p\in M$, and $A_{\eta}$ be the second
fundamental tensor with respect to $\eta$.
\begin{defn}
The manifold $M$ is said to be h-umbilical with respect to $\eta$ if
$$A_{\eta}=\frac{{\rm trace}\, A_{\eta}}{2n}\, I-
\frac{{\rm trace}\, (A_{\eta}\circ J)}{2n}\, J.$$
\end{defn}

Using Lemma \ref {L: 2.1}, it follows immediately
\begin{lem}\label {L: 2.4}
If at a point $p\in M$ \,$\sigma = g\,H- \tilde g\, JH$, then $M$ is h-umbilical with respect
to every vector $\eta$ normal to $M$ at $p$.
\end{lem}

\begin{lem}\label {L: 2.5}
If there exists a vector $\eta$ normal to $M$ at $p\in M$ such that $M$ is h-umbilical with
respect to $\eta$, then $\sigma = g\, H- \tilde g\, JH$ at $p$.
\end{lem}

If $M$ is umbilical, i.e. $\sigma = g\, H$, and $M$ is in the same time h-umbilical, then
$H=0$ and $M$ is totally geodesic.

From Lemma \ref {L: 2.1} and the equalities (2.7), (2.8), (2.9) we obtain
\begin{lem}\label {L: 2.6}
Let $M$ be h-umbilical. Then
\begin{itemize}
\item[1)] if $M$ is totally real, i.e. $H\neq 0,\; H\perp JH$, then $M$ is umbilical with
respect to $\xi$ or $J\xi$, i.e. $A_{\xi}=\displaystyle{\frac{1}{2n}\,({\rm trace}\,A)\, I}$ or
$A_{J\xi}=\displaystyle{\frac{1}{2n}\,({\rm trace}\,(A\circ J))\, I}$;
\vskip 2mm
\item[2)] if $H$ is isotropic, i.e. $g(H, H)=0$, then $M$ is umbilical with respect to $H$ or
$JH$, i.e. $\displaystyle{A_H=\frac{1}{n}({\rm trace}\,A)\, I}$ or
$\displaystyle{A_{JH}=\frac{1}{n}({\rm trace}\,(A\circ J))\, I}$.
\end{itemize}
\end{lem}

The h-umbilical holomorphic hypersurfaces have the following property.
\begin{thm}\label {T: 2.2}
Let $M'\; (\dim M'= 2n+2\geq 8)$ be of constant totally real sectional curvatures $\nu\,'$
and $\tilde\nu\,'$. If $M$ is h-umbilical, then $M$ is of constant totally real sectional
curvatures $\nu$ and $\tilde \nu$ such that
$$\nu=\nu\,' + g(H, H),\quad \tilde\nu = \tilde\nu\,' + \tilde g(H, H).$$
\end{thm}

{\it Proof.} Since $M'$ is of constant totally real sectional curvatures $\nu\,'$ and
$\tilde\nu\,'$, Theorem \ref {T: 2.1} implies
$$R'=\nu\,'\,(\pi_1-\pi_2) + \tilde\nu\,'\,\pi_3.$$
Using Lemma \ref {L: 2.4}, we find
$A=\displaystyle{\frac{1}{2n}\,({\rm trace}\,A)\, I -
\frac{1}{2n}\,({\rm trace}\,(A\circ J))\, J}$. Substituting $R'$ and $A$ into the first
equality of (2.4), we obtain
$$R=\left\{\nu\,'+\left(\frac{{\rm trace}\,A}{2n}\right)^2-
\left(\frac{{\rm trace}\,(A\circ J)}{2n}\right)^2\right\}\,(\pi_1-\pi_2)
+\left\{\tilde\nu\,' +2\;\frac{{\rm trace\,A}}{2n}\;
\frac{{\rm trace\,(A\circ J)}}{2n}\right\}\,\pi_3.$$
Applying Theorem \ref {T: 2.1} we get that $M$ is of constant totally real sectional
curvatures $\nu$ and $\tilde\nu$. Taking into account (2.8), we find
$$\nu=\nu\,' + g(H, H),\quad \tilde\nu = \tilde\nu\,' + \tilde g(H, H).$$
\hfill{QED}

\begin{cor}\label {C:2.1}
Let $M'$ \,$(\dim M'=2n+2\geq 8)$ be of constant totally real sectional curvatures. If $M$ is
connected and h-umbilical, then on $M$
\begin{itemize}
\item[(1)] $g(H, H)= const,\quad \tilde g(H, H)= const;$
\vskip 2mm
\item[(2)] ${\rm trace}\, A = const,\quad {\rm trace}\,(A\circ J) = const.$
\end{itemize}
\end{cor}

\begin{cor}\label {C: 2.2}
Let $M'$\, $(\dim M'=2n+2\geq 6)$ be of constant totally real sectional curvatures $\nu\,'$
and $\tilde\nu\,'$. If $M$ is totally geodesic, then $M$ is of constant totally real sectional
curvatures $\nu\,'$ and $\tilde\nu\,'$.
\end{cor}

\vskip 2mm
\section{Examples of K\"ahler manifolds with Norden metric of constant totally real sectional
curvatures}

Let $M'=\mathbb{R}^{2n+2}$ be equipped with the canonical complex structure $J$ and the metric
$g$, given by (1.1). Then $(M', g, J)$ is a K\"ahler manifold with Norden metric. In this case
the curvature tensor $R'$ of $M'$ is zero. Identifying the point
$z=(x^1,..., x^{n+1}; y^1,...,y^{n+1})$ in $M'$ with the position vector $Z$, we define the
submanifold $S^{2n}(z_0; a, b)$ by the equalities
$$g(Z-Z_0, Z-Z_0)=a,\leqno (3.1)$$
$$\tilde g(Z-Z_0, Z-Z_0)=b,\leqno (3.2)$$
where $a, b \in \mathbb{R},\; (a, b)\neq (0, 0)$.

$S^{2n}$ is a $2n$-dimensional submanifold
of $M'$ and the vectors $Z-Z_0$ and $ J(Z-Z_0)$ are perpendicular to $T_z S^{2n}$. The condition
$(a, b)\neq (0, 0)$ implies that the rank of $g$ on $T_z S^{2n}$ is equal to $2n$ and
$T_z S^{2n}$ is $J$-invariant, i.e. $S^{2n}$ is a holomorphic hypersurface of $M'$.

\begin{defn}
The holomorphic hypersurface $S^{2n}(z_0; a, b)$ of $\mathbb{R}^{2n}$ determined by the
equalities (3.1) and (3.2) is said to be an h-sphere with a center $z_0$ and parameters
$a$, $b$.
\end{defn}

In the case $a=1,\, b=0$\; $S^{2n}(z_0; 1, 0)$ is the sphere of Kotel'nikov-Study \cite {3}.
\vskip 2mm

Let us consider the vector fields $\{\xi, J\xi\}$ normal to $S^{2n}$ and satisfying (2.1).
We have
$$\xi=-\lambda Z- \mu JZ,\quad J\xi = \mu Z- \lambda JZ.\leqno (3.3)$$
Taking into account (3.1) and (3.2) we obtain
$$\lambda^2-\mu^2=\frac{a}{a^2+b^2},\quad 2 \lambda \mu= \frac{b}{a^2 +b^2}.\leqno (3.4)$$

Thus, for $\lambda$ and $\mu$ we have
$$\lambda^2=\frac{\sqrt{a^2+b^2}+ a}{2(a^2+b^2)},\quad
\mu^2=\frac{\sqrt{a^2+b^2}- a}{2(a^2+b^2)}, \quad sign\,(\lambda \mu) = sign\,  b.$$

Let $X$ be a vector in $T_z S^{2n}$. Using (3.2) we find
$$-A X= \nabla'_X \xi = -\lambda \nabla'_X Z - \mu J\nabla'_X Z.$$
Since $\nabla'$ is flat, then $\nabla'_X Z = X$, $Z$ being the position vector field. Hence,
$A=\lambda I+ \mu J$. The last formula implies that $S^{2n}$ is h-umbilical with respect to
$\xi$ and
$$\lambda = \frac{1}{2n}\,{\rm trace}\,A,\quad \mu= -\frac{1}{2n}\,{\rm trace}\,(A\circ J).$$
Applying Lemma \ref {L: 2.5}, Theorem \ref {T: 2.2} and (3.4), we obtain

\begin{thm}\label {T: 3.1}
Every h-sphere $S^{2n}(z_0; a, b)$ in $\mathbb{R}^{2n+2}\, (2n+2 \geq 6)$ is an h-umbilical
holomorphic hypersurface of constant totally real sectional curvatures
$$\nu=\frac{a}{a^2+b^2}= g(H, H),\quad \tilde\nu=\frac{-b}{a^2+b^2}=\tilde g(H, H).
\leqno (3.5)$$
\end{thm}

Thus, we obtain holomorphic hypersurfaces of $\mathbb{R}^{2n+2}$ with prescribed constant
totally real sectional curvatures.

\begin{thm}\label {T: 3.2}
 Let $\nu$ and $\tilde \nu$ \, $(\nu^2+\tilde\nu^2 >0)$ be real numbers. Then
 $\displaystyle{S^{2n}\left(z_0; \frac{\nu}{\nu^2+\tilde\nu^2},
 \frac{-\tilde\nu}{\nu^2+\tilde\nu^2}\right)}$\; $(n\geq 2)$ is an h-sphere in
 $\mathbb{R}^{2n+2}$ of constant totally real sectional curvatures $\nu$ and $\tilde\nu$.
\end{thm}
\vskip 2mm

The following special cases are worth to be noted:
\begin{itemize}
\item[1)] $S^{2n}(z_0; a, 0)$.\;
In this case $\nu = (1/ a),\; \tilde\nu = 0$. Hence, the sphere
of Kotel'nikov-Study has constant totally real sectional curvatures $\nu=1,\, \tilde\nu=0$.
\vskip 2mm
\item[2)] $S^{2n}(z_0; 0, b)$.\;
In this case $\nu = 0,\; \tilde\nu = -(1/ b)$.
\end{itemize}
\vskip 2mm

\begin{defn}
The h-sphere $\bar S^{2n}(z_0; a, -b)$ is said to be conjugate to the h-sphere
$S^{2n}(z_0; a, b)$.
\end{defn}

If $S^{2n}$ has constant totally real sectional curvatures $\nu$ and $\tilde\nu$, than
Theorem \ref {T: 3.2} and (3.5) imply $\bar S^{2n}$ has constant totally real sectional
curvatures $\nu$ and $-\tilde\nu$.

If $M$ is a holomorphic hyperplane in $\mathbb{R}^{2n+2}$, then $M$ is of totally real
sectional curvatures $\nu=\tilde\nu=0$.

The hypersurface of $\mathbb{R}^{2n+2}$ given by the equations (3.1) and (3.2) with
$a=b=0$ is an example of a holomorphic isotropic hypersurface of $\mathbb{R}^{2n+2}$.

\vskip 2mm
\section{A classification of the holomorphic hypersurfaces of constant totally real
sectional curvatures in $\mathbb{R}^{2n+2}$}

In the previous section we proved that the holomorphic hyperplanes and the h-spheres have
constant totally real sectional curvatures. In this section we consider the inverse question.
\begin{thm}\label {T: 4.1}
Let $M'$  $(\dim M'= 2n+2\geq 8)$ be a K\"ahler manifold with Norden metric of constant
totally real sectional curvatures $\nu\,'$ and $\tilde\nu\,'$. If $M$ is a holomorphic
hypersurface of $M'$ with constant totally real sectional curvatures $\nu$ and $\tilde\nu$,
then $M$ is h-umbilical.
\end{thm}

{\it Proof.} Since $M'$ and $M$ are of constant totally real sectional curvatures Theorem
\ref {T: 2.1} and Lemma \ref {L: 2.3} imply
$$2 A^2 y = ({\rm trace}\,A)\, Ay-({\rm trace}\,(A\circ J))\, JAy +
2(n-1)\{(\nu\,'-\nu)\, y+(\tilde\nu\,'-\tilde\nu)\,Jy\} \leqno (4.1)$$
for arbitrary $y\in T_p M,\,p\in M$.

Applying Lemma \ref {L: 1.1} to the h-symmetric operator $A$, we get
$$Ax_k=\lambda_k\,x_k+\mu_k\,Jx_k,\quad AJx_k=-\mu_k\,x_k+\lambda_k\,Jx_k,
\; k=1,..., n,\leqno (4.2)$$
where $\{x_1,..., x_n; Jx_1,..., Jx_n\}$ is an adapted basis for $T_p M,\, p\in M$, of
h-proper vectors of $A$.

Using the first formula of (2.4) and (4.2) we find
$$\nu\,'-\nu=\mu_j\,\mu_k-\lambda_j\,\lambda_k,\quad
\tilde\nu\,'-\tilde\nu=\lambda_j\,\mu_k+\lambda_k\,\mu_j,\; j\neq k.\leqno (4.3)$$
Substituting (4.2) in (4.1) and taking into account (4.3) we check
$$\begin{array}{l}
{\rm trace}\,A=2(\lambda_j+(n-1)\lambda_k)=2(\lambda_k+(n-1)\lambda_j);\; j\neq k,\\
[2mm]
{\rm trace}\,(A\circ J)=-2(\mu_j+(n-1)\mu_k)=-2(\mu_k+(n-1)\mu_j);\; j\neq k.
\end{array}$$

Since $n\geq 3$ these equalities imply $\lambda_j=\lambda,\; \mu_j=\mu$ for all $j=1,..., n$.
Hence, because of (4.2)
$$A=\displaystyle{\frac{1}{2n}\,({\rm trace}\,A)\, I -
\frac{1}{2n}\,({\rm trace}\,(A\circ J))\, J}.$$

Applying  Lemma \ref {L: 2.5} we obtain the proposition.
\hfill{QED}
\vskip 2mm

\begin{cor}\label {C: 4.1}
Let $M'$ and $M$ be as in Theorem $\ref {T: 4.1}$. $M$ is totally geodesic if and only if
$\nu\,'=\nu$ and $\tilde\nu\,'=\tilde\nu$.
\end{cor}

The statement follows from (4.3) taking into account $\lambda_j=\lambda$, $\mu_j=\mu$
$(j=1,..., n)$.
\begin{cor}\label {C: 4.2}
Let $M$ be a holomorphic hypersurface of $\mathbb{R}^{2n+2}$ $(2n+2\geq 8)$ with constant
totally real sectional curvatures $\nu$ and $\tilde \nu$. Then $M$ is h-umbilical and
$g(H, H)=\nu,\; \tilde g(H, H)=\tilde\nu$.
\end{cor}

\begin{cor}\label {C: 4.3}
Let $M$ be a holomorphic hypersurface of $\mathbb{R}^{2n+2}$ $(2n+2\geq 8)$ with zero
totally real sectional curvatures. Then $M$ is totally geodesic.
\end{cor}

\begin{thm}\label {T: 4.2}
Let $M$ be a connected holomorphic h-umbilical hypersurface in
$\mathbb{R}^{2n+2}$ $(2n+2\geq 8)$. Then $M$ lies on an h-sphere or on a holomorphic
hyperplane.
\end{thm}

{\it Proof.} Let $\mathbb{U}$ be a coordinate neighborhood on $M$ and $\{\xi, J\xi\}$
be normal vector fields on $\mathbb{U}$, satisfying (2.1). From the condition of the
theorem we have $A=\lambda I+\mu J$ on $\mathbb{U}$, where
$$\lambda = \frac{1}{2n}\,{\rm trace}\,A,\quad \mu= -\frac{1}{2n}\,{\rm trace}\,(A\circ J).$$
Theorem \ref {T: 2.2} implies that $M$ is of constant totally real sectional curvatures
$\nu= \lambda^2-\mu^2,$ $\tilde\nu=-2\lambda \mu$. Hence $\lambda= const,$ $\mu= const$ on
$\mathbb{U}$.

Identifying $z\in \mathbb{U}$ with the position vector $Z$ in $\mathbb{R}^{2n+2}$ we
consider the vector field $\xi+\lambda Z+\mu JZ$ on $\mathbb{U}$. Using the Weingarten
formula and taking into account that $\nabla'$ is flat, we obtain
$$\nabla'_x (\xi+\lambda Z+\mu JZ)=-AX+\lambda X +\mu JX=0$$
for arbitrary vector field $X$ on $\mathbb{U}$ tangent to $M$. The last equality implies
$$\xi+\lambda Z+\mu JZ = const.\leqno (4.4)$$

Let $\lambda^2+\mu^2>0$. Then (4.4) can be written in the form
$$\xi+\lambda Z+\mu JZ=\lambda Z_0+\mu JZ_0, \quad Z_0= const.$$
From here we find
$$Z-Z_0=\frac{-\lambda \xi+\mu J\xi}{\lambda^2+\mu^2},\quad
J(Z-Z_0)=\frac{-\mu \xi-\lambda J\xi}{\lambda^2+\mu^2}.$$
Thus, for every $z\in \mathbb{U}$
$$g(Z-Z_0, Z-Z_0)=\frac{\lambda^2-\mu^2}{(\lambda^2+\mu^2)^2},\quad
\tilde g(Z-Z_0, Z-Z_0)=\frac{2\lambda \mu}{(\lambda^2+\mu^2)^2}.$$
Hence, $\mathbb{U}$ lies on a h-sphere $S^{2n}$ with a center $z_0$ and parameters
$\displaystyle{\frac{\lambda^2-\mu^2}{(\lambda^2+\mu^2)^2}}$,\;
$\displaystyle{\frac{2\lambda \mu}{(\lambda^2+\mu^2)^2}}$.

So we proved that for every $z\in M$ there exists $\mathbb{U}\ni z$ such that $\mathbb{U}$
lies on an h-sphere $S^{2n}$.

Let $\mathbb{U}$ be a fixed coordinate neighborhood and $S^{2n}$ be the h-sphere such that
$\mathbb{U}\subset S^{2n}$. $M_0$ will stand for the set of points $z$ in $M$ belonging to
$S^{2n}$ together with a neighborhood $\mathbb{V}_z\ni z$. Obviously $M_0\neq \emptyset$
and $M_0$ is open.

Let $w$ be in the closure $ \overline{M_0}$ of $M_0$. Then there exists $\mathbb{W}\ni w$ and
$\mathbb{W}$ lies on an h-sphere $S_1^{2n}$. The open set $\mathbb{W}\cap M_0\neq \emptyset$
lies on $S_1^{2n}$ and on $S^{2n}$. Hence, $S_1^{2n}=S^{2n}$ and $M_0=M$ because of the
connectedness of $M$.

In the case $\lambda=\mu=0$ from (4.4) it follows $\xi= const$ on $\mathbb{U}$. If $X$ is
a vector field on $\mathbb{U}$ tangent to $M$ and $z\in \mathbb{U}$, then
$\nabla'_X g(Z, \xi)=g(X, \xi)=0$. From here it follows that
$$g(Z, \xi)= const=g(Z_0, \xi),\quad \tilde g(Z, \xi)=\tilde g(Z_0, \xi)$$
and $\mathbb{U}$ lies on the holomorphic hyperplane
$$g(\xi, Z-Z_0)=0,\quad \tilde g(\xi, Z-Z_0)=0.\leqno (4.5)$$
Further, as in the previous case, it follows that $M$ lies on the holomorphic hyperplane (4.5).
\hfill{QED}
\vskip 2mm

Using theorems \ref {T: 4.1} and \ref {T: 4.2} we obtain the following classification theorem.
\begin{thm}
Let $M$ be a connected holomorphic hypersurface of $(\mathbb{R}^{2n+2}, g, J)\, (n+1\geq 4)$.
If $M$ is of constant totally real sectional curvatures $\nu$ and $\tilde \nu$, then
\begin{itemize}
\item[1)] $M$ lies on an h-sphere $S^{2n}$ with parameters
$\displaystyle{ \frac{\nu}{\nu^2+\tilde\nu^2},\;  \frac{-\tilde\nu}{\nu^2+\tilde\nu^2}}$,
when $\nu^2+\tilde\nu^2>0$;
\vskip 2mm
\item[2)] $M$ lies on a holomorphic hyperplane, when $\nu=\tilde\nu=0$.
\end{itemize}
\end{thm}

\end{document}